\documentclass{amsproc}

\newtheorem{theorem}{Theorem}[section]
\newtheorem{lemma}[theorem]{Lemma}

\newtheorem{conjecture}{Conjecture}[section]

\theoremstyle{definition}
\newtheorem{definition}[theorem]{Definition}
\newtheorem{example}[theorem]{Example}

\theoremstyle{remark}
\newtheorem{remark}[theorem]{Remark}

\numberwithin{equation}{section}

\begin{document}

\title{Derived equivalence of holomorphic symplectic manifolds}

\author{Justin Sawon}
\address{Department of Mathematics, SUNY at Stony Brook, New York
11794} 
\email{sawon@math.sunysb.edu}
\thanks{The author was supported in part by NSF Grant \#0305865.}

\date{October 31, 2003.}



\begin{abstract}
We use twisted Fourier-Mukai transforms to study the relation between
an abelian fibration on a holomorphic symplectic manifold and its dual
fibration. Our reasoning leads to an equivalence between the derived
category of coherent sheaves on one space and the derived category
of twisted sheaves on the other space.
\end{abstract}

\maketitle

\section{Introduction}

Fourier-Mukai transforms were introduced by Mukai~\cite{mukai81} in
the early 80s, and have recently been instrumental in the
understanding of several important problems in higher dimensional
algebraic geometry (for example, see Bridgeland, King, and
Reid~\cite{bkr01} or Bridgeland~\cite{bridgeland02}). The original
example of Mukai related the derived categories of an abelian variety
and its dual. In this paper we investigate a relative version of this
equivalence using the twisted Fourier-Mukai transforms of C{\u a}ld{\u
a}raru~\cite{caldararu00}; this leads to (twisted) derived
equivalences for holomorphic symplectic manifolds.

The main idea is to begin with a holomorphic symplectic manifold which
is fibred by abelian varieties (which we call an {\em abelian
fibration\/}), construct its dual fibration, and then relate the
derived category of coherent sheaves on the original manifold to the
derived category of twisted sheaves on the dual fibration. The
twisting comes from a gerbe which arises as the natural obstruction to
extending the fibrewise equivalence to a family. 

The motivation for studying this problem came from a desire to
understand the relation between different abelian fibrations. An
elliptic K3 surface can be deformed, through elliptic surfaces, to an
elliptic K3 surface with a section. For a higher dimensional abelian
fibration $X$ we'd like to construct a fibration $X^0$ which is
locally isomorphic to $X$ and which admits a section (so that $X$ is a
torsor over $X^0$). We'd also like $X^0$ to be a deformation of $X$,
though the example in Subsection 5.4 indicates that this may not
always be the case. Elliptic K3 surfaces which admit sections are
easily classified via their Weierstra{\ss} models, and the hope is
that in higher dimensions one may be able to classify (holomorphic
symplectic) abelian fibrations which admits sections (see
Sawon~\cite{sawon03}).

Section 2 contains some preliminaries on holomorphic gerbes, twisted
sheaves, and derived categories. In Section 3 we review Fourier-Mukai
transforms, their twisted version, and some results important for
their application. In Section 4 we discuss fibrations by abelian
varieties on holomorphic symplectic manifolds. We then describe how
one can construct a dual fibration and a derived equivalence coming
from a twisted Fourier-Mukai transform. Singular fibres provide a
significant technical obstacle: even in the `nicest' cases our results
depend on several as-yet unverified assumptions (these are explained
in the text; one involves extending an autoduality result to
compactified Jacobians, the other involves extending a Fourier-Mukai
result to the twisted case). In the final section we present examples
of abelian fibrations and their duals. The author intends to present
the first of these examples in complete rigour in a future
article~\cite{sawon??}. The remaining examples are meant to illustrate
some interesting behaviour that warrants further investigation.

The author would like to thank the organizers of the {\em Workshop on
algebraic structures and moduli spaces\/} for putting together such an 
interesting schedule of talks, and for support to attend the workshop
and present this work. Tom Bridgeland, Andrei C{\u a}ld{\u a}raru,
Eduardo de Sequeira Esteves, Manfred Lehn, Daisuke Matsushita, and
Richard Thomas have also been very helpful in answering the author's
questions.

\section{Preliminaries}

\subsection{Holomorphic gerbes}

In trying to extend certain local sheaves to global objects we will
come up against obstructions in the form of holomorphic gerbes. Gerbes
are higher dimensional analogues of line bundles, as we now
explain. Let $M$ be a complex manifold with Stein cover $\{U_i\}$. An
element in $\mathrm{H}^1(M,\mathcal{O}^*)$ is an equivalence class of
holomorphic line bundles on $M$. A holomorphic gerbe is the geometric
object whose equivalence classes are elements of
$\mathrm{H}^2(M,\mathcal{O}^*)$.

\begin{definition}
A (holomorphic) gerbe is a collection $\{\mathcal{L}_{ij}\}$ of line
bundles on two-fold intersections $U_{ij}:=U_i\cap U_j$ such that
\begin{enumerate}
\item $\mathcal{L}_{ii}=\mathcal{O}_{U_i},$
\item $\mathcal{L}_{ji}=\mathcal{L}_{ij}^{-1},$
\item
$\mathcal{L}_{ijk}:=\mathcal{L}_{ij}\otimes\mathcal{L}_{jk}\otimes\mathcal{L}_{ki}$
has a given trivialization,
\item the trivialization of
$\mathcal{L}_{ijk}\otimes\mathcal{L}_{jkl}^{-1}\otimes\mathcal{L}_{kli}\otimes\mathcal{L}_{lij}^{-1}$
induced by the given trivializations of $\mathcal{L}_{ijk}$ is the
canonical one.
\end{enumerate}
\end{definition}

\begin{remark}
Since we chose our cover to be Stein, each line bundle
$\mathcal{L}_{ij}$ is trivializable. If we fix trivializations of
each $\mathcal{L}_{ij}$, this gives a trivialization of
$\mathcal{L}_{ijk}$, but this may differ by a non-zero holomorphic
function 
$$\beta_{ijk}\in\Gamma(U_{ijk},\mathcal{O}^*_{U_{ijk}})$$
from the trivialization given in the definition. The fourth condition
then states that $\delta\beta$ is trivial, so that $\beta$ is a
cocycle representing a class in $\mathrm{H}^2(M,\mathcal{O}^*)$.
\end{remark}

\begin{definition}
The isomorphism class of a gerbe is the class
$[\beta]\in\mathrm{H}^2(M,\mathcal{O}^*)$. We will usually abuse
notation by writing simply $\beta$ for $[\beta]$.
\end{definition}

\begin{remark}
Isomorphism of gerbes can be defined in the expected way, by taking
into account different choices of trivialization of
$\mathcal{L}_{ij}$ and refinement of the cover (for a more detailed
discussion see Hitchin~\cite{hitchin01}). There is a bijection between
the set of all gerbes isomorphic to a given gerbe and the set of all
cocycles representing $[\beta]$.
\end{remark}

In this article, all our constructions can be made to depend only on
the isomorphism class of the gerbe involved, not the gerbe
itself. We will often abuse terminology by referring to the
isomorphism class simply as a gerbe.

\begin{example}
The exponential exact sequence
$$0\rightarrow\mathbb{Z}\rightarrow\mathcal{O}\rightarrow\mathcal{O}^*\rightarrow
0$$
induces the long exact sequence
$$\ldots\rightarrow\mathrm{H}^2(M,\mathbb{Z})\rightarrow\mathrm{H}^2(M,\mathcal{O})\rightarrow\mathrm{H}^2(M,\mathcal{O}^*)\rightarrow\mathrm{H}^3(M,\mathbb{Z})\rightarrow\ldots$$

On a Calabi-Yau three-fold $\mathrm{H}^2(M,\mathcal{O})$ vanishes and
therefore $\mathrm{H}^2(M,\mathcal{O}^*)$ injects into the lattice
$\mathrm{H}^3(M,\mathbb{Z})$ and is discrete. We say that each
topological gerbe has a unique holomorphic structure.

On a K3 surface $\mathrm{H}^2(M,\mathcal{O})=\mathbb{C}$ and
$\mathrm{H}^3(M,\mathbb{Z})$ vanishes. Therefore
$\mathrm{H}^2(M,\mathcal{O}^*)$ is a quotient of $\mathbb{C}$
(possibly non-Hausdorff, depending on the rank of the image of
$\mathrm{H}^2(M,\mathbb{Z})$). We say that there are no non-trivial
topological gerbes, only non-trivial holomorphic structures.
\end{example}

As explained by C{\u a}ld{\u a}raru~\cite{caldararu00}, a gerbe can
also be described as a sheaf of Azumaya algebras over $M$. However,
this applies only to gerbes whose class in
$\mathrm{H}^2(M,\mathcal{O}^*)$ is torsion, a restriction that we'd
prefer to avoid.

\subsection{Twisted sheaves}

Suppose we are given a gerbe, denoted $\beta$, on the complex manifold
$M$ with corresponding trivializations
$$\beta_{ijk}\in\Gamma(U_{ijk},\mathcal{O}^*_{U_{ijk}}).$$
There is a natural notion of sheaves twisted by $\beta$.

\begin{definition}
A $\beta$-twisted sheaf on $M$, denoted $\mathcal{F}$, is a collection
$\{\mathcal{F}_i\}$ of local coherent sheaves (on $U_i$) and
isomorphisms $\phi_{ij}:\mathcal{F}_i\rightarrow\mathcal{F}_j$
satisfying
\begin{enumerate}
\item $\phi_{ji}=\phi_{ij}^{-1}$,
\item $\phi_{ki}\circ\phi_{jk}\circ\phi_{ij}=\beta_{ijk}\mathrm{Id}$.
\end{enumerate}
\end{definition}

\begin{remark}
The second condition means that the composition on the left, which is
an automorphism of $\mathcal{F}_i$ restricted to $U_{ijk}$, is given
by $\beta_{ijk}$ times the identity isomorphism.
\end{remark}

If there exists a locally free $\beta$-twisted sheaf of rank $n$, then
one can show that $\beta$ must be torsion of order dividing $n$. In
this case, when $\beta$ is torsion, twisted sheaves can also be
describes as sheaves of modules over the corresponding sheaf of
Azumaya algebras (see C{\u a}ld{\u a}raru~\cite{caldararu00}).

The space of gerbes forms a group which we write multiplicatively. We
can also pull back a gerbe $\beta$ by a map $f:N\rightarrow M$ to get
a gerbe $f^*\beta$ on $N$.

\begin{lemma}[C{\u a}ld{\u a}raru~\cite{caldararu00}]
The pull-back of a $\beta$-twisted sheaf $\mathcal{F}$ on $M$ by
$f:N\rightarrow M$ is a $f^*\beta$-twisted sheaf $f^*\mathcal{F}$ on
$N$. The tensor product of a $\beta_1$-twisted sheaf $\mathcal{F}_1$
and a $\beta_2$-twisted sheaf $\mathcal{F}_2$ on $M$ is a
$\beta_1\otimes\beta_2$-twisted sheaf
$\mathcal{F}_1\otimes\mathcal{F}_2$ on $M$. In particular, if
$\beta_2=\beta_1^{-1}$, the inverse gerbe, then
$\mathcal{F}_1\otimes\mathcal{F}_2$ is untwisted, i.e.\ a genuine
sheaf on $M$.
\end{lemma}

Strictly speaking, the tensor product makes sense for locally free
sheaves; for more general sheaves we should pass to the derived
category and use derived tensor product, which we shall do in the next
subsection.

\subsection{Derived categories}

The basic idea behind working with derived categories is that an
object in an abelian category should be identified with all its
resolutions. This ensures that certain functors behave nicely. For
example, exactness of sequences of coherent sheaves is not preserved
under tensor product or push-forward, but if we restrict to projective
or injective sheaves (respectively) then exactness is preserved. In
this case, rather than working with the abelian category of coherent
sheaves, we should take its derived category so that we can replace
sheaves by projective or injective resolutions. See
Thomas~\cite{thomas01} or page 143 of Gelfand and Manin~\cite{gm96}
for further motivation behind the derived category.

\begin{definition}
Let $\mathcal{C}$ be an abelian category (basically, the set of
morphisms between any two objects is an abelian group, and kernels and 
cokernels exist; see Definition 4.1 in~\cite{thomas01}). Examples
include the category $\mathrm{Coh}(M)$ of coherent sheaves on a
complex manifold $M$, or in the presence of a gerbe $\beta$, the
category $\mathrm{Coh}(M,\beta)$ of $\beta$-twisted sheaves.

Let $\mathrm{Kom}(\mathcal{C})$ be the category of chain complexes
over $\mathcal{C}$. A quasi-isomorphism is a morphism between two
complexes which induces an isomorphism on the level of cohomology. 

The derived category $\mathcal{D}(\mathcal{C})$ of $\mathcal{C}$ is
the category whose objects are chain complexes over $\mathcal{C}$ and
whose morphisms are usual morphisms of chain complexes, plus formal
inverses of quasi-isomorphisms. Thus in the derived category,
quasi-isomorphic objects become isomorphic. The bounded derived
category is given by considering only bounded chain complexes.

In this article we will use the bounded derived categories of
$\mathrm{Coh}(M)$ and $\mathrm{Coh}(M,\beta)$ on $M$, which we denote
by $\mathcal{D}^b_{\mathrm{coh}}(M)$ and
$\mathcal{D}^b_{\mathrm{coh}}(M,\beta)$ respectively.
\end{definition}

The derived category is an additive category, but not abelian as
kernels and cokernels do not exist. Instead, given a morphism
$B^{\bullet}\rightarrow C^{\bullet}$ of chain complexes, there exists
an ``exact triangle''
$$A^{\bullet}\rightarrow B^{\bullet}\rightarrow C^{\bullet}\rightarrow
A^{\bullet}[1]$$
where $A^{\bullet}[1]$ is the complex $A^{\bullet}$ shifted one place
to the right. It is these exact triangles which are preserved under
tensor product and push-forward, and which induce long exact sequences
in cohomology. Thus derived categories are examples of triangulated
categories, and when we talk about functors between derived categories
we shall always mean functors which preserve exact triangles.

The derived category $\mathcal{D}^b_{\mathrm{coh}}(M)$ of a projective
variety $M$ has recently emerged as an important tool in studying the
birational geometry of $M$. It contains a great deal of information
about $M$. Indeed Bondal and Orlov~\cite{bo01} proved that if $M$ has
ample canonical or anti-canonical bundle, it can be recovered up to
isomorphism from its derived category. This is no longer the case when
the canonical bundle is trivial. In the case of K3 surfaces,
Orlov~\cite{orlov97} showed that two K3 surfaces have equivalent
derived categories if and only if there is a Hodge isometry between
their transcendental lattices.

In this article we deal with holomorphic symplectic manifolds, higher
dimensional analogues of K3 surfaces. Our aim is to find pairs of
non-isomorphic manifolds with equivalent derived categories; more
accurately, the derived category of one manifold should be equivalent
to the twisted derived category of the other (so this is a not quite a
direct analogue of the K3 case, as studied by Orlov).

\section{Fourier-Mukai transforms}

\subsection{Untwisted}

Fourier-Mukai transforms were first introduced in~\cite{mukai81} as a
way of producing equivalences of derived categories of non-isomorphic
varieties. Recently there have been some striking applications of
these methods, including Bridgeland, King, and Reid's ``derived'' McKay
correspondence~\cite{bkr01} and Bridgeland's construction of
three-fold flops~\cite{bridgeland02}.

Let $M$ be a fine moduli space of stable sheaves on some projective
variety $X$. In other words, $M$ parametrizes a complete family of
stable sheaves on $X$ and there exists a universal sheaf $\mathcal{U}$
on $X\times M$. Denote by $\pi_X$ and $\pi_M$ the projections to $X$
and $M$ respectively.

\begin{definition}[Mukai~\cite{mukai81}]
We call the functor
$$\Phi^{\mathcal{U}}_{M\rightarrow
X}:\mathcal{D}^b_{\mathrm{coh}}(M)\rightarrow\mathcal{D}^b_{\mathrm{coh}}(X)$$
$$\mathcal{E}^{\bullet}\mapsto\mathbf{R}\pi_{X*}(\mathcal{U}\stackrel{\mathbf{L}}{\otimes}\pi_M^*\mathcal{E}^{\bullet})$$
an integral transform. If it is an equivalence of categories we call
it a Fourier-Mukai transform.
\end{definition}

\begin{remark}
Defining this operation on the category of coherent sheaves would be
unsatisfactory, as it would not preserve exactness. Instead we use
derived categories, and derived functors which preserve exact
triangles.
\end{remark}

\begin{remark}
The dual sheaf $\mathcal{U}^{\vee}$ allows us to define a functor
$$\Phi^{\mathcal{U}^{\vee}}_{X\rightarrow
M}:\mathcal{D}^b_{\mathrm{coh}}(X)\rightarrow\mathcal{D}^b_{\mathrm{coh}}(M),$$
which is the inverse equivalence in the case of a Fourier-Mukai
transform. Note that the dimensions of $X$ and $M$ must agree if their
derived categories are equivalent.
\end{remark}

\begin{example}
Let $X$ be an elliptic curve $E$, and let $M$ the dual elliptic curve 
$\widehat{E}$, which parametrizes degree zero line bundles on $E$. The
Poincar{\'e} line bundle $\mathcal{P}$ is a universal bundle on
$E\times\widehat{E}$. In higher dimensions, we could let $X$ be an
abelian variety and $M$ its dual. Once again the Poincar{\'e} line
bundle is a universal bundle. Mukai~\cite{mukai81} showed these
universal bundles induce equivalences of derived categories.
\end{example}

Let $\mathcal{O}_m$ be the skyscraper sheaf supported at $m\in
M$. Then $\mathcal{U}_m:=\Phi^{\mathcal{U}}_{M\rightarrow
X}\mathcal{O}_m$ is the sheaf on $X$ which the point $m\in M$
parametrizes. The skyscraper sheaves on $M$ `span' the derived
category and in some sense are `orthonormal' with respect to the
$\mathrm{Ext}^{\bullet}$-pairing. If the integral transform
$\Phi^{\mathcal{U}}_{M\rightarrow X}$ is an equivalence it should
preserve these properties. Based on these ideas, the following theorem
was proved by Bridgeland, extending work of Mukai, Bondal, and Orlov.

\begin{theorem}[Bridgeland~\cite{bridgeland99}]
\label{bridgeland_criteria}
Assume that $M$ has the same dimension as $X$, and both are smooth
varieties. The functor $\Phi^{\mathcal{U}}_{M\rightarrow X}$ is an
equivalence of categories if and only if the following criteria are
satisfied:
\begin{enumerate}
\item for all $m\in M$, $\mathcal{U}_m$ is a simple sheaf, i.e.\
$$\mathrm{Hom}(\mathcal{U}_m,\mathcal{U}_m)=\mathbb{C},$$
and $\mathcal{U}_m\otimes\mathcal{K}_X=\mathcal{U}_m$ where
$\mathcal{K}_X$ is the canonical bundle of $X$,
\item for all $m_1\neq m_2\in M$ and all integers $i$,
$$\mathrm{Ext}^i_X(\mathcal{U}_{m_1},\mathcal{U}_{m_2})=0.$$
\end{enumerate}
\end{theorem}

\begin{example}
In the previous example, when $X$ is an abelian variety and $M$ its
dual, $\mathcal{U}_m$ is a line bundle and
$$\mathrm{Ext}^i_X(\mathcal{U}_{m_1},\mathcal{U}_{m_2})=\mathrm{H}^i(X,\mathcal{U}_{m_1}^{\vee}\otimes\mathcal{U}_{m_2})$$
vanishes for $m_1\neq m_2$ by standard results on the cohomology of
line bundles on abelian varieties (see chapter 3 of Birkenhake and
Lange~\cite{bl92}). The first criterion is also easy to verify, and
hence the integral transform is indeed an equivalence.
\end{example}

In general the second criterion can be more difficult to
verify. However, Bridgeland and Maciocia proved the following theorem,
which says it is enough to check the condition on a set of sufficiently
large codimension.

\begin{theorem}[Bridgeland-Maciocia~\cite{bm02}]
\label{bridgeland_maciocia}
Assume that $M$ has the same dimension $n$ as $X$, and that $X$ is a 
smooth variety. Suppose that the following criteria are satisfied: 
\begin{enumerate}
\item for all $m\in M$, $\mathcal{U}_m$ is a simple sheaf and
$\mathcal{U}_m\otimes\mathcal{K}_X=\mathcal{U}_m$,
\item for all $m_1\neq m_2\in M$
$$\mathrm{Hom}_X(\mathcal{U}_{m_1},\mathcal{U}_{m_2})=0,$$
and the closed subscheme
$$\Gamma(\mathcal{U}):=\{(m_1,m_2)\in M\times
M|\mathrm{Ext}^i_X(\mathcal{U}_{m_1},\mathcal{U}_{m_2})\neq 0\mbox{
for some }i\in\mathbb{Z}\}$$
of $M\times M$ has dimension at most $n+1$.
\end{enumerate}
Then in fact $\Gamma(\mathcal{U})$ is the diagonal, $M$ is also
smooth, and $\Phi^{\mathcal{U}}_{M\rightarrow X}$ is an equivalence of
categories. 
\end{theorem}

\begin{remark}
If the dimension $n$ is greater than two, it can also be difficult to
show directly that $M$ is smooth; but the theorem takes care of that.
\end{remark}

\subsection{Twisted}

In~\cite{caldararu00} C{\u a}ld{\u a}raru generalized the ideas of the
previous subsection to twisted sheaves by introducing twisted
Fourier-Mukai transforms. Before giving the definition, let us explain
how twisted sheaves naturally arise in this context.

Perhaps the clearest motivation comes from considering relative moduli
spaces. Suppose that $p_X:X\rightarrow B$ is a proper fibration and
let $M$ be an irreducible component of the relative moduli space of
stable sheaves on $X$. In other words, each point of $M$ corresponds
to a stable sheaf supported on a fibre of $X\rightarrow B$. Thus $M$
is also fibred over $B$, and we write $p_M:M\rightarrow B$.

\begin{example}
Let $X\rightarrow B$ be an elliptic fibration and let $M$ be its
compactified relative Jacobian $\overline{\mathrm{Jac}}(X/B)$ (see
D'Souza~\cite{dsouza79}, Altman-Iarrobino-Kleiman~\cite{aik77}, or
Rego~\cite{rego80}). Then $M$ parametrizes torsion-free rank one
degree zero sheaves on the fibres of $X$. Since these sheaves are rank
one, a destabilizing sheaf cannot exist, and hence $M$ is a component
of the moduli space of stable sheaves on $X$ (we use Simpson's
definition of stability and construction of the moduli
space~\cite{simpson94}).

To construct $M$ in this example, smooth fibres of $X$ are
replaced by their dual elliptic curves. Of course elliptic curves are
self-dual, but the global structure of $M$ may also differ from
$X$. For instance, $M$ has a canonical section (given by the trivial
sheaf on each fibre of $X$) whereas $X$ may not admit any global
section.
\end{example}

Let $t\in B$, and consider the fibres $X_t$ and $M_t$. One can show
that provided $X_t$ is reduced (it need not be irreducible) then $M_t$
is a moduli space of stable sheaves on $X_t$. This may not be the case
if $X_t$ is non-reduced because $M_t$ might include `fat' sheaves
whose scheme-theoretic support is non-reduced. We assume that $M_t$ is
a fine moduli space of sheaves on $X_t$, so that a universal sheaf
exists, and moreover that these fibrewise universal sheaves fit
together locally to give a local universal sheaf over an open subset
in $B$ containing $t$.

\begin{example}
In the previous example, the Poincar{\'e} line bundle is a universal
sheaf for any smooth fibre of $X_t$ and its dual $M_t$. One can extend
it over a neighbourhood $B_i$ of $t\in B$.
\end{example}

Now assume that $B$ is covered by open sets $B_i$ such that there
exists a local universal sheaf $\mathcal{U}_i$ over each $X\times
M_i$, where $M_i:=p_M^{-1}(B_i)$. More precisely, the local universal
sheaf is, a priori, really a sheaf over $X_i\times_{B_i}M_i$; but we
can extend it (by zero) to $X\times_B M_i$ and then push it forward by
the inclusion of the fibre-product in $X\times M_i$. On the
intersection $X\times M_{ij}$ of $X\times M_i$ and $X\times M_j$, the
restrictions of the two sheaves $\mathcal{U}_i$ and $\mathcal{U}_j$
are universal sheaves for the same moduli problem. Hence they are
isomorphic and there exists a line bundle $\mathcal{L}_{ij}$ on
$M_{ij}$ such that
$$\mathcal{U}_i|_{X\times
M_{ij}}=\pi_M^*\mathcal{L}_{ij}\otimes\mathcal{U}_j|_{X\times
M_{ij}}.$$
The fibre of $\mathcal{L}_{ij}$ over a point $m\in M_{ij}$ is
identified with the set of morphisms from $\mathcal{U}_i|_{X\times m}$
to $\mathcal{U}_j|_{X\times m}$. These sheaves are both isomorphic to
$\mathcal{U}_m$, and thus we have implicitly used the fact that $M$
parametrizes stable, and hence simple sheaves on $X$; i.e.\ the fibres
of $\mathcal{L}_{ij}$ are isomorphic to $\mathbb{C}$.

\begin{theorem}[C{\u a}ld{\u a}raru~\cite{caldararu00}]
\label{obstruction}
In the situation described above, the collection $\{\mathcal{L}_{ij}\}$
defines a gerbe on $M$. The class
$\beta\in\mathrm{H}^2(M,\mathcal{O}^*)$ of this gerbe represents the
obstruction to patching the local universal sheaves $\mathcal{U}_i$
together; i.e.\ there exists a global universal sheaf on $X\times M$
if and only if $\beta=0$.
\end{theorem}

\begin{example}
\label{elliptic_fibration}
For an elliptic fibration $X\rightarrow B$ and its compactified
relative Jacobian $M:=\overline{\mathrm{Jac}}(X/B)$ one can show that
the following are equivalent (provided the singular fibres are well
enough behaved):
\begin{enumerate}
\item $X$ is isomorphic to $M$,
\item $X$ admits a global section,
\item there is a global universal sheaf on $X\times M$,
\item $\beta$ vanishes.
\end{enumerate}
The first two statements are equivalent because $X$ is a torsor over
$M$. The last two statements are equivalent by the previous
theorem. For the remaining equivalence, one shows that a local
section of $X\rightarrow B$ can be used to locally extend the
Poincar{\'e} bundle, thereby producing a local universal bundle (see
C{\u a}ld{\u a}raru~\cite{caldararu00} for details).
\end{example}

Next observe that the collection $\{\mathcal{U}_i\}$ of local
universal sheaves forms a $\pi_M^*\beta$-twisted sheaf on $X\times M$,
which we denote simply by $\mathcal{U}$. So in the twisted case, the
failure of the existence of a universal sheaf is precisely controlled
by $\beta$, and moreover we have instead a {\em twisted universal
sheaf\/} $\mathcal{U}$. We can also construct integral transforms with
twisted sheaves.

\begin{definition}
If the functor
$$\Phi^{\mathcal{U}}_{M\rightarrow
X}:\mathcal{D}^b_{\mathrm{coh}}(M,\beta^{-1})\rightarrow\mathcal{D}^b_{\mathrm{coh}}(X)$$
$$\mathcal{E}^{\bullet}\mapsto\mathbf{R}\pi_{X*}(\mathcal{U}\stackrel{\mathbf{L}}{\otimes}\pi_M^*\mathcal{E}^{\bullet})$$
is an equivalence of categories we call it a twisted Fourier-Mukai
transform.
\end{definition}

\begin{remark}
In the definition, $\mathcal{E}^{\bullet}$ denotes a complex of
$\beta^{-1}$-twisted sheaves on $M$. Therefore
$\pi_M^*\mathcal{E}^{\bullet}$ is a complex of
$\pi_M^*\beta^{-1}$-twisted sheaves on $X\times M$, and
$\mathcal{U}\stackrel{\mathbf{L}}{\otimes}\pi_M^*\mathcal{E}^{\bullet}$
is an untwisted sheaf on $X\times M$. The final push-forward takes us
to the derived category of $X$.
\end{remark}

C{\u a}ld{\u a}raru~\cite{caldararu00} generalized Bridgeland's
criteria for when an integral transform is an equivalence
(Theorem~\ref{bridgeland_criteria}) to the twisted case. Skyscraper
sheaves can always be regarded as twisted sheaves, so defining
$\mathcal{U}_m:=\Phi^{\mathcal{U}}_{M\rightarrow X}\mathcal{O}_m$
still makes perfect sense in the twisted case. The criteria are then
precisely the same as in the untwisted case, and we will not repeat
them here.

\begin{example}
\label{elliptic_equivalence}
If $\beta$ is non-zero in Example~\ref{elliptic_fibration}, then
provided the singular fibres are well enough behaved, we can apply the
criteria and hence show that
$$\Phi^{\mathcal{U}}_{M\rightarrow X}:\mathcal{D}^b_{\mathrm{coh}}(M,\beta^{-1})\rightarrow\mathcal{D}^b_{\mathrm{coh}}(X)$$
is an equivalence. Various examples where $X$ is an elliptic Calabi-Yau
three-fold are discussed in C{\u a}ld{\u
a}raru~\cite{caldararu00}. Dealing with the singular fibres is one of
the biggest difficulties, so usually one begins with an example whose
singular fibres are fairly mild.
\end{example}

\section{Holomorphic symplectic manifolds}

\subsection{Abelian fibrations}

\begin{definition}
Let $X$ be a compact K{\"a}hler manifold. We call $X$ a holomorphic
symplectic manifold if it admits a closed non-degenerate holomorphic
two-form
$$\sigma\in\mathrm{H}^0(X,\Lambda^2T^*)=\mathrm{H}^{2,0}(X)$$
($\sigma$ is the holomorphic analogue of a symplectic structure). Note
that the canonical bundle $\mathcal{K}_X$ is trivialized by
$\sigma^{\wedge n}$, where $n$ is half the dimension of $X$. If $X$ is
simply connected and $\mathrm{H}^{2,0}(X)=\mathbb{C}$ is generated by
$\sigma$ then we say $X$ is irreducible.
\end{definition}

\begin{remark}
By Bogomolov's decomposition theorem a holomorphic symplectic manifold
has a finite unramified cover which is a product of complex tori and
irreducible holomorphic symplectic manifolds. The latter are therefore
the {\em building blocks\/} of the theory, and we will always mean
`irreducible holomorphic symplectic' when we write simply `holomorphic
symplectic', unless stated otherwise.
\end{remark}

In two dimensions the only (irreducible) holomorphic symplectic
manifolds are K3 surfaces. A proper morphism from a K3 surface to a
curve is necessarily an elliptic fibration, with base $\mathbb{P}^1$
and generic fibre an elliptic curve. In higher dimensions we have the
following theorem of Matsushita.

\begin{theorem}[Matsushita~\cite{matsushita99,matsushita00}]
Let $X$ be a projective irreducible holomorphic symplectic manifold of
dimension $2n$. Let $p_X:X\rightarrow B$ be a proper surjective
projective morphism, with connected generic fibre, and with projective
base $B$ of dimension strictly between zero and $2n$. Then
\begin{enumerate}
\item the generic fibre is a (holomorphic) Lagrangian abelian variety
of dimension $n$,
\item the base is Fano with the same Hodge numbers as $\mathbb{P}^n$.
\end{enumerate}
It follows from the second statement and the classification of
surfaces that the base must be $\mathbb{P}^2$ when $n=2$.
\end{theorem}

\begin{definition}
If there exists a morphism $p_X:X\rightarrow B$ as in the theorem, we
call $X$ an abelian fibration. In this case $X$ is also known as an
algebraically complete integrable system. 
\end{definition}

A K3 surface is elliptic if and only if it contains a non-trivial
divisor with square zero (with respect to intersection
pairing). Similarly, in higher dimensions we expect that $X$ is an
abelian fibration if and only if it contains a divisor with some
special properties. Presently this is only conjectural (see
Sawon~\cite{sawon03}). Moreover, we don't yet know whether every
holomorphic symplectic manifold can be deformed to an abelian
fibration.

\subsection{Dual fibrations}

Now we come to the heart of this article. Given an abelian fibration,
our aim is to construct its dual fibration and then relate it, via a
twisted Fourier-Mukai transform, to the original fibration. In this
subsection we outline the general theory, which is a higher
dimensional analogue of the elliptic fibrations discussed in
Examples~\ref{elliptic_fibration} and~\ref{elliptic_equivalence}. In
the next section we will present some examples that indicate what kind
of behaviour we can expect.

Let $p_X:X\rightarrow B$ be an abelian fibration of dimension $2n$. We
will assume that all the fibres are geometrically integral (i.e.\
reduced and irreducible).

\begin{definition}
The compactified relative Picard scheme
$\overline{\mathrm{Pic}}^0(X/B)$ of $X$ is the moduli space
parametrizing torsion-free rank one sheaves of degree zero (vanishing
first Chern class) supported on the fibres of $X$. Its construction,
which is due to Altman and Kleiman~\cite{ak80}, already requires that
the fibres of $X$ be geometrically integral.
\end{definition}

Let $M$ be the compactified relative Picard scheme
$\overline{\mathrm{Pic}}^0(X/B)$ of $X$. Then $M$ parametrizes stable
sheaves: they are pure-dimension with geometrically integral support,
so there can be no destabilizing subsheaf (we are using Simpson's
definition of stability~\cite{simpson94}). 

A priori we don't know very much about $M$. It need not be a smooth
space. There is a projection $p_M:M\rightarrow B$, though we cannot
expect this to be flat in general, as we now explain. If $X_t$ is a
smooth fibre of $X$, $t\in B$, then clearly the corresponding fibre
$M_t$ of $M$ is simply the dual abelian variety. More generally, if
$X_t$ is a singular fibre one can show that $M_t$ parametrizes sheaves
on $X$ of the form $\iota_*\mathcal{E}$, where
$\iota:X_t\hookrightarrow X$ is the inclusion of the fibre and
$\mathcal{E}$ is a stable rank one degree zero sheaf on $X_t$. In
particular, $M_t$ is always isomorphic to a moduli space of stable
sheaves on $X_t$ (our singular fibres are geometrically integral,
though in fact this appears to be true provided $X_t$ contains no
non-reduced components).

This gives us a precise description of the singular fibres of $M$,
but unfortunately moduli spaces of stable sheaves on singular
varieties have not been extensively studied. Undesirable behaviour
cannot automatically be ruled out: for example, the singular fibres
of $M$ need not be irreducible and they may contain components of
dimension greater than $n$. One example we can try to understand is
the (compactified) Jacobian of a degeneration of a curve.

\begin{example}
\label{autoduality}
Let $C$ be a singular curve contained in a surface. Assume that $C$ 
has at worst double point singularities (for example, nodes, cusps, or
tacnodes). Altman, Kleiman, and Iarrobino~\cite{aik77} showed that the
compactified Jacobian $J:=\overline{\mathrm{Jac}}C$ of such a curve is
geometrically integral; it is a degeneration of an abelian
variety. We can then apply Altman and Kleiman's construction to $J$,
and hence obtain the compactified Picard scheme
$\overline{\mathrm{Pic}}^0J$. A recent theorem of Esteves, Gagn{\'e},
and Kleiman~\cite{egk02} gives an isomorphism of uncompactified Picard
schemes
$$\mathrm{Pic}^0C\cong\mathrm{Pic}^0J.$$
It is expected~\cite{esteves03} that this isomorphism will extend to
the compactifications
$$\overline{\mathrm{Pic}}^0C\cong\overline{\mathrm{Pic}}^0J.$$
The left hand side is just $J$ itself, up to isomorphism.

In terms of our abelian fibration $X$ and its dual fibration $M$, this
means that if the singular fibre $X_t$ can be identified as the
compactified Jacobian of a mildly singular curve then the
corresponding singular fibre $M_t$, which is
$\overline{\mathrm{Pic}}^0(X_t)$, will be isomorphic to $X_t$. In
particular $M_t$ is geometrically integral and has dimension $n$, so
$M\rightarrow B$ will be flat. Singular fibres like this will arise in
the examples of the next section.
\end{example}

The example shows that in some special cases we can acquire an exact
understanding of the singular fibres of $M$, so we will proceed under
the assumption that the singular fibres of $M$ are sufficiently
well-behaved. In fact, we'd like $M$ to be locally isomorphic to $X$
as a fibration over $B$. However, this is too much to ask for in
general: if the fibres of $X$ are not principally polarized then even
a smooth fibre $X_t$ of $X$ need not be isomorphic to the
corresponding fibre $M_t$ of $M$. This problem does not arise with
elliptic fibres, which are always principally polarized. The remedy is
to define $p_0:X^0\rightarrow B$ to be the compactified relative
Picard scheme $\overline{\mathrm{Pic}}^0(M/B)$ of $p_M:M\rightarrow
B$. For a smooth fibre $X_t$ of $X$, the corresponding fibre $X^0_t$
of $X^0$ is the double dual of, and hence canonically isomorphic to,
$X_t$. We will assume the same is true for singular fibres (this
should be true when $X_t$ is a mild degeneration of a Jacobian, as in
the last example).

Thus $p_X:X\rightarrow B$ and $p_0:X^0\rightarrow B$ are locally
isomorphic as fibrations. Since the compactified relative Picard
scheme over $B$ always has a canonical section, $X$ is therefore a
torsor over $X^0$.

\begin{lemma}
The set of torsors $X$ over $X^0$ is one-to-one correspondence with
the sheaf cohomology $\mathrm{H}^1(B,X^0)$, where by $X^0$ we really
mean the sheaf of local sections of $X^0\rightarrow B$. 
\end{lemma}

\begin{proof}
Choose an open cover (in the analytic topology) $\{B_i\}$ of $B$
such that $X$ and $X^0$ are isomorphic fibrations over each $B_i$,
i.e.\ there exists local isomorphisms
$$f_i:X_i:=p_X^{-1}(B_i)\rightarrow X^0_i:=p_0^{-1}(B_i)$$
over $B_i$. On the overlap $X^0_{ij}:=p_0^{-1}(B_{ij})$ the
composition $f_j\circ f_i^{-1}$ gives an automorphism $\alpha_{ij}$ of
the fibration $X^0_{ij}\rightarrow B_{ij}$, which is given by a
translation in each fibre. Since we have a basepoint in each fibre
($X^0$ has a global section), the family of translations $\alpha_{ij}$
is equivalent to a section of $X^0_{ij}\rightarrow B_{ij}$. These
local sections form a 1-cocyle
$\alpha\in\mathrm{H}^1(B,X^0)$. Conversely, choosing a representative
of $\alpha$ gives us the families of translations $\alpha_{ij}$ on
overlaps needed to glue the local fibrations $X^0_i\rightarrow B_i$
together to obtain $X$.
\end{proof}

\begin{lemma}
To each torsor $X$ over $X^0$ we can assign a gerbe
$\beta\in\mathrm{H}^2(M,\mathcal{O}^*)$ on $M$. This is the same gerbe
that arises as the obstruction to the existence of a global universal
sheaf on $X\times M$.
\end{lemma}

\begin{proof}
We have seen how to construct $\alpha\in\mathrm{H}^1(B,X^0)$ from
$X$. Taking a representative, we obtain local sections $\alpha_{ij}$
of $X^0$; but each point of $X^0$ represents a rank one degree zero
sheaf (generically a line bundle) supported on a fibre of
$M\rightarrow B$, so the local section $\alpha_{ij}$ gives a line
bundle $\mathcal{L}_{ij}$ on $M_{ij}:=p_M^{-1}(B_{ij})$. One can show
that this collection of line bundles forms a gerbe
$\beta\in\mathrm{H}^2(M,\mathcal{O}^*)$, which does not depend on our
choice of representative of $\alpha$.

To relate this to the gerbe arising as an obstruction in
Theorem~\ref{obstruction}, we first consider $X^0$ as a moduli space
on $M$. On an abelian variety and its dual, the Poincar{\'e} line
bundle is a universal sheaf. Since both $X^0$ and $M$ have global
sections, these Poincar{\'e} line bundles fit together to give a
global universal sheaf on $M\times X^0$. Restricting to 
$$M\times X^0_i\cong M\times X_i$$
gives a bundle whose support is contained in $M_i\times
X_i$. Interchanging the factors, and extending by zero, gives a local
universal sheaf on $X\times M_i$. Note that we have implicitly used
the relation between the Poincar{\'e} bundle of an abelian variety and
the Poincar{\'e} bundle of its dual (see Exercise 16 on page 45 of
Birkenhake and Lange~\cite{bl92}). Finally, it is not difficult to see
that these local universal sheaves no longer patch together to give a
global sheaf on $X\times M$; moreover they differ on overlaps
precisely by tensoring with $\mathcal{L}_{ij}$.
\end{proof}

\begin{remark}
The lemma implies that there is a map 
$$\mathrm{H}^1(B,X^0)\rightarrow\mathrm{H}^2(M,\mathcal{O}^*).$$
Our next theorem says that this is an inclusion, but it will not be an
isomorphism in general because the gerbes we obtain from torsors $X$
over $X^0$ have a somewhat specialized form.
\end{remark}

Now we come to our main result, the generalization of
Example~\ref{elliptic_fibration} to abelian fibrations.

\begin{theorem}
\label{abelian_tfae}
Let $X\rightarrow B$ be an abelian fibration, $M$ the compactified
relative Picard scheme of $X$, and $X^0$ the compactified relative
Picard scheme of $M$. Assume that for all $t\in B$ the fibres $X^0_t$
and $X_t$ are isomorphic. Let $\beta\in\mathrm{H}^2(M,\mathcal{O}^*)$
be the gerbe associated to the torsor $X$ over $X^0$, as in the
previous lemma. Then the following are equivalent:
\begin{enumerate}
\item $X$ is isomorphic to $X^0$,
\item $X$ admits a global section,
\item there is a global universal sheaf on $X\times M$,
\item $\beta$ vanishes.
\end{enumerate}
\end{theorem}

\begin{proof}
The first two statements are equivalent because $X$ is a torsor over
$X^0$. The last two statements are equivalent by the previous lemma
and Theorem~\ref{obstruction}. The remaining equivalence is shown
using the same method as in the elliptic case: namely one shows that a
local section of $X\rightarrow B$ can be used to locally extend the
Poincar{\'e} bundle, thereby producing a local universal bundle.
\end{proof}

Finally we want to generalize Example~\ref{elliptic_equivalence}, and
show that we have an equivalence of derived categories for abelian
fibrations. So assume we are in the same situation as in
Theorem~\ref{abelian_tfae}, and that $\beta$ is non-zero. We want to
use C{\u a}ld{\u a}raru's generalization of Bridgeland's criteria to
show that the twisted universal sheaf $\mathcal{U}$ induces an
equivalence
$$\Phi^{\mathcal{U}}_{M\rightarrow X}:\mathcal{D}^b_{\mathrm{coh}}(M,\beta^{-1})\rightarrow\mathcal{D}^b_{\mathrm{coh}}(X).$$
For all $m\in M$, $\mathcal{U}_m$ is a stable and hence simple sheaf
on $X$. Moreover $X$ is holomorphic symplectic, with trivial canonical
bundle, and thus
$\mathcal{U}_m\otimes\mathcal{K}_X=\mathcal{U}_m$. The first criterion
is therefore satisfied, and it remains to show that if $m_1\neq m_2\in
M$ then
$$\mathrm{Ext}^i_X(\mathcal{U}_{m_1},\mathcal{U}_{m_2})=0$$
for all integers $i$. We can consider various cases:
\begin{enumerate}
\item if $m_1$ and $m_2$ belong to different fibres of $M$, then
$\mathcal{U}_{m_1}$ and $\mathcal{U}_{m_2}$ are supported on different
fibres of $X$, and the $\mathrm{Ext}$-groups must vanish,
\item if $m_1$ and $m_2$ belong to the same smooth fibre of $M$, then
$\mathcal{U}_{m_1}$ and $\mathcal{U}_{m_2}$ are (push-forwards of)
distinct line bundles $L_1$ and $L_2$ on the same smooth fibre $X_t$
of $X$; in this case all cohomology groups
$\mathrm{H}^i(X_t,L_1^*\otimes L_2)$ vanish (see Chapter 3 of
Birkenhake and Lange~\cite{bl92}), the normal bundle to the fibre is
holomorphically trivial, and the Koszul spectral sequence then shows
that the $\mathrm{Ext}$-groups also vanish (as in Subsection 7.2 of
Bridgeland and Maciocia~\cite{bm02}),
\item if $m_1$ and $m_2$ belong to the same singular fibre of $M$, it
is more difficult to proceed with this calculation.
\end{enumerate}
For elliptic fibrations, a direct approach to the third case is
feasible. For abelian fibrations, we instead want to apply Bridgeland
and Maciocia's Theorem~\ref{bridgeland_maciocia}. Since the proof of
that theorem uses only local arguments, the result should still be
valid in the twisted case. We need to bound the dimension of
$\Gamma(\mathcal{U})$. Let $X$ have dimension $2n$. Singular fibres
occur over a codimension one subset of the base $B$, and the fibres
themselves have dimension $n$, so the pairs $(m_1,m_2)$ corresponding
to case three above are contained in a subset of $M\times M$ of
dimension $3n-1$. The theorem requires this to be at most one larger
than the dimension of $X$, and therefore $n$ can be at most two.

Using these arguments, we can expect to produce equivalences
$$\Phi^{\mathcal{U}}_{M\rightarrow X}:\mathcal{D}^b_{\mathrm{coh}}(M,\beta^{-1})\rightarrow\mathcal{D}^b_{\mathrm{coh}}(X).$$
for holomorphic symplectic four-folds $X$ fibred by abelian surfaces,
under some additional hypotheses. In higher dimensions the same
equivalence ought to hold, but there are some additional technical
difficulties yet to be overcome. We will discuss some examples of
abelian fibrations in the next section.

\section{Examples}

\subsection{The Hilbert scheme of two points on a K3}

\begin{definition}
Let $S$ be a K3 surface. The Hilbert scheme of two points on $S$ is
the space
$$S^{[2]}:=\mathrm{Blow}_{\Delta}(S\times S)/{\mathbb{Z}_2}$$
obtained by blowing up the diagonal in $S\times S$ and quotienting by
the involution interchanging the two factors.
\end{definition}

Fujiki~\cite{fujiki83} showed that the four-fold $S^{[2]}$ has a
holomorphic symplectic structure; this was the first example of a
higher dimensional irreducible holomorphic symplectic manifold. We
will be interested in certain deformations of $S^{[2]}$ which are
abelian fibrations, i.e.\ fibred over $\mathbb{P}^2$ by abelian
surfaces.

\begin{example}
Let $S$ be a double cover of the plane $(\mathbb{P}^2)^{\vee}$
ramified over a generic sextic. The linear system of lines in
$(\mathbb{P}^2)^{\vee}$ is simply the dual plane
$\mathbb{P}^2$. Pulling these lines back to $S$ gives a complete
linear system of genus two curves on $S$. Denote the family of these
curves by $\mathcal{C}\rightarrow\mathbb{P}^2$. One can show that all
of these curves are irreducible (they develop at worst two nodes or a
cusp), and we can therefore consider the compactified relative
Jacobian $Z^2:=\overline{\mathrm{Jac}}(\mathcal{C}/\mathbb{P}^2)$, and
more generally the compactified degree $d$ relative Picard scheme
$Z^d:=\overline{\mathrm{Pic}}^d(\mathcal{C}/\mathbb{P}^2)$. Each space
$Z^d$ is clearly a fibration
$$\begin{array}{ccc}
\overline{\mathrm{Pic}}^d & \hookrightarrow & Z^d \\
 & & \downarrow \\
 & & \mathbb{P}^2 \\
\end{array}$$
with generic fibre a smooth abelian surface.

Now given a degree $d$ line bundle (or more generally, a torsion-free
rank one sheaf) on a curve $C$ in the linear system $\mathbb{P}^2$, we
can take the push-forward by the inclusion $C\hookrightarrow S$. This
gives a (torsion) sheaf on $S$, which moreover is stable. We can
therefore regard $Z^d$ as an irreducible component of the moduli space
of stable sheaves on $S$. Mukai~\cite{mukai84} showed that these
moduli spaces are smooth and admit holomorphic symplectic forms, and
hence $Z^d$ is an irreducible holomorphic symplectic four-fold and an
abelian fibration. Moreover, it is a deformation of the Hilbert scheme
$S^{[2]}$ (see Yoshioka~\cite{yoshioka99}, for instance).
\end{example}

We summarize some further properties of the four-folds $Z^d$, which
hold for all $d\in\mathbb{Z}$ (these statements are mostly
straight-forward, and will be proved in Sawon~\cite{sawon??}): 
\begin{itemize}
\item all fibres of $Z^d$ are irreducible,
\item $Z^d$ is locally isomorphic to $Z^{d+1}$ as a fibration over
$\mathbb{P}^2$,
\item $Z^0\rightarrow\mathbb{P}^2$ admits a global section, given by
taking the trivial degree zero line bundle on each curve,
\item $Z^1\rightarrow\mathbb{P}^2$ admits a birational two-valued
section,
\item $Z^d$ is globally isomorphic to $Z^{d+2}$, the isomorphism given
by tensoring each degree $d$ line bundle by the canonical bundle of
the genus two curve; thus $Z^d$ admits a global section if $d$ is
even, and a birational two-valued section if $d$ is odd,
\item $Z^1$ is not isomorphic to $Z^0$, as they have different
weight-two Hodge structures; thus $Z^d$ does not admit a section if
$d$ is odd.
\end{itemize}

Thus we essentially have just two different spaces, with $Z^1$ a
torsor over $Z^0$, and $Z^1$ corresponds to a two-torsion element
$\alpha\in\mathrm{H}^1(\mathbb{P}^2,Z^0)$. These spaces were
investigated by Markushevich~\cite{markushevich95,markushevich96}; he
showed that a holomorphic symplectic four-fold which is also the
compactified relative Jacobian of a family of curves must be
isomorphic to $Z^0$. Our approach will be slightly different, as we
next want to interpret $Z^0$ as a moduli space of sheaves on the {\em
four-fold\/} $Z^1$, so that we can set up a twisted Fourier-Mukai
transform between them.

\begin{theorem}
\label{compact_Picard}
Under a mild assumption (explained in the proof below), $Z^0$ is the
compactified relative Picard scheme
$\overline{\mathrm{Pic}}^0(Z^1/\mathbb{P}^2)$ of
$Z^1\rightarrow\mathbb{P}^2$.
\end{theorem}

\begin{proof}
The fibres of $Z^1$ are compactified Jacobians of surficial curves
with at worst double point singularities. As explained in
Example~\ref{autoduality}, we expect that the autoduality result of
Esteves, Gagn{\'e}, and Kleiman~\cite{egk02} should extend to
compactified Picard schemes, and thus the fibres of
$\overline{\mathrm{Pic}}^0(Z^1/\mathbb{P}^2)$ will be isomorphic to
the corresponding fibres of $Z^1$, and of $Z^0$. Thus
$\overline{\mathrm{Pic}}^0(Z^1/\mathbb{P}^2)$ will be a torsor over
$Z^0$; since they both admit global sections they must be
isomorphic.
\end{proof}

\begin{remark}
Thus $Z^0$ is the dual fibration of $Z^1$. We also obtain the
smoothness of the compactified relative Picard scheme of $Z^1$. This
is difficult to prove directly, as the compactified relative Picard
scheme is a moduli space of sheaves on a four-fold, and moduli spaces
on higher dimensional varieties are poorly understood.
\end{remark}

There exists a gerbe $\beta\in\mathrm{H}^2(Z^0,\mathcal{O}^*)$ and a
twisted universal sheaf $\mathcal{U}$ on $Z^1\times Z^0$. Since $Z^0$
and $Z^1$ are not isomorphic, $\beta$ is non-trivial by
Theorem~\ref{abelian_tfae}.

\begin{theorem}
\label{equivalence}
Assuming that the theorem of Bridgeland and Maciocia
(Theorem~\ref{bridgeland_maciocia}) extends to the twisted case, we
obtain an equivalence of derived categories 
$$\Phi^{\mathcal{U}}_{Z^0\rightarrow Z^1}:\mathcal{D}^b_{\mathrm{coh}}(Z^0,\beta^{-1})\rightarrow\mathcal{D}^b_{\mathrm{coh}}(Z^1).$$
\end{theorem}

\begin{proof}
This follows from the arguments at the end of Subsection 4.2. In this
example the degeneracy locus is a curve in $\mathbb{P}^2$, so pairs of
points $(m_1,m_2)$ in the same singular fibre contribute a five
dimensional component to $\Gamma(\mathcal{U})\subset Z^0\times
Z^0$. Since $\mathrm{dim}Z^1+1=5$, $\Gamma(\mathcal{U})$ is not too
large.
\end{proof}

\subsection{The Hilbert scheme of more than two points on a K3}

Generalizing the example of the last subsection, the Hilbert scheme
$S^{[n]}$ of $n$ points on a K3 surface $S$ is a minimal
desingularization of the symmetric product
$\mathrm{Sym}^nS$. Beauville~\cite{beauville83} showed that $S^{[n]}$
is an irreducible holomorphic symplectic manifold of dimension $2n$.

\begin{example}
Let $S$ be a genus $g$ K3 surface, $g>2$, i.e.\ $S$ contains a smooth
genus $g$ curve $C$. Assume that $S$ is otherwise generic. The curve
moves in a $g$ dimensional linear system; denote this family of curves
by $\mathcal{C}\rightarrow\mathbb{P}^g$. One can show that all the
curves are irreducible, and hence we can consider the compactified
degree $d$ relative Picard scheme
$Z^d:=\overline{\mathrm{Pic}}^d(\mathcal{C}/\mathbb{P}^g)$. These give
fibrations
$$\begin{array}{ccc}
\overline{\mathrm{Pic}}^d & \hookrightarrow & Z^d \\
 & & \downarrow \\
 & & \mathbb{P}^g \\
\end{array}$$
whose generic fibres are $g$ dimensional abelian varieties.

As in the previous subsection, the spaces $Z^d$ can also be
interpreted as irreducible components of the Mukai moduli space of
stable sheaves on $S$~\cite{mukai84}. Thus they are smooth irreducible
holomorphic symplectic manifolds. They have dimension $2g$ and are
deformation equivalent to the Hilbert scheme $S^{[g]}$
(see~\cite{yoshioka99}), which can therefore be deformed to an abelian
fibration.
\end{example}

We can generalize some of the statements of the previous subsection:
\begin{itemize}
\item all fibres of $Z^d$ are irreducible,
\item $Z^d$ is locally isomorphic to $Z^{d+1}$ as a fibration over
$\mathbb{P}^g$,
\item $Z^0\rightarrow\mathbb{P}^g$ admits a global section, given by
taking the trivial degree zero line bundle on each curve,
\item $Z^d$ is globally isomorphic to $Z^{d+2g-2}$, the isomorphism
given by tensoring each degree $d$ line bundle by the canonical bundle
of the genus $g$ curve; thus $Z^d$ admits a global section if $d$ is a
multiple of $2g-2$.
\end{itemize}

Thus $Z^1$, $Z^2$, $\ldots$, and $Z^{2g-1}$ are all torsors over
$Z^0$, and there are corresponding elements $\alpha_1$, $\alpha_2$,
$\ldots$, and $\alpha_{2g-1}\in\mathrm{H}^1(\mathbb{P}^g,Z^0)$. These $2g-1$
spaces won't all be different; for example, there are isomorphisms
$$Z^{2g-1}\cong Z^{-1}\cong Z^1$$
where the second isomorphism is given by taking a line bundle of
degree $-1$ to its dual. Nevertheless, it should be possible to show
that we obtain at least some spaces which are not isomorphic to
$Z^0$ (for instance, by looking at their weight-two Hodge
structures). Let $Z^i$ be such a space.

\begin{conjecture}
The space $Z^0$ is the compactified relative Picard scheme
$\overline{\mathrm{Pic}}^0(Z^i/\mathbb{P}^g)$ of
$Z^i\rightarrow\mathbb{P}^g$. Hence there is a (non-trivial) gerbe
$\beta\in\mathrm{H}^2(Z^0,\mathcal{O}^*)$ and a twisted universal
sheaf $\mathcal{U}$ on $Z^i\times Z^0$, which induces an equivalence
of derived categories 
$$\Phi^{\mathcal{U}}_{Z^0\rightarrow Z^i}:\mathcal{D}^b_{\mathrm{coh}}(Z^0,\beta^{-1})\rightarrow\mathcal{D}^b_{\mathrm{coh}}(Z^i).$$
\end{conjecture}

\begin{remark}
We first need to look carefully at the curves in the family
$\mathcal{C}\rightarrow\mathbb{P}^g$ to see what kind of singularities
they can acquire; the first statement should then follow from a
compact version of the autoduality result of Esteves, Gagn{\'e}, and
Kleiman~\cite{egk02}, as in the proof of Theorem~\ref{compact_Picard}.

To show that we get an equivalence we need the twisted version of
Bridgeland and Maciocia's theorem, as in the proof of
Theorem~\ref{equivalence}. However, this is not enough since the
degeneracy locus is codimension one in $\mathbb{P}^g$, so the
dimension of $\Gamma(\mathcal{U})$ is
$$(g-1)+2\mathrm{dim}\{\mbox{fibre of }Z^0\}=3g-1>2g+1=\mathrm{dim}Z^i+1$$
as $g>2$. To make $\Gamma(\mathcal{U})$ smaller, we need to show
directly that
$$\mathrm{Ext}^j_{Z^i}(\mathcal{U}_{m_1},\mathcal{U}_{m_2})=0$$
for all integers $j$ in some of the cases when $m_1$ and $m_2$ belong
to the same singular fibre of $Z^0$. This can be done: the generic
point $m_1$ on a singular fibre of $Z^0$ parametrizes a locally free
sheaf on the corresponding singular fibre of $Z^i$, and then some
perseverance with the standard machinery of homological algebra shows
that the above $\mathrm{Ext}$-groups vanish (the same can be done when
$m_2$ parametrizes a locally free sheaf). This leaves the case when
$m_1$ and $m_2$ both parametrize non-locally free sheaves, so both
$m_1$ and $m_2$ are codimension one in the fibre of $Z^0$. This gives
$$\mathrm{dim}\Gamma(\mathcal{U})=3g-3\leq
2g+1=\mathrm{dim}Z^i+1$$
provided $g\leq 4$. We expect the equivalence should still exist when
$g>4$, but new techniques must be developed to verify this. 
\end{remark}

\subsection{The ten dimensional example of O'Grady}

In the previous subsection we considered a K3 surface containing a
genus $g$ curve, which is otherwise generic. The K3 surface $S$ of
Section 5.1 (the double cover of the plane) contains a genus five
curve $E$, the pull-back of a conic from the plane; however, it is not
otherwise generic as it also contains a genus two curve. So we don't
expect the conjecture of the previous subsection to
apply. Nevertheless, this example is still of interest as we now
explain.

The curve $E$ moves in a five dimensional linear system
$\mathbb{P}^5$; we denote the family of curves by
$\mathcal{E}\rightarrow\mathbb{P}^5$. These curves are certainly not
all irreducible: for example, the family includes the pull-backs from
the plane of singular conics (pairs of lines). The family also
includes non-reduced curves (pull-backs of double lines). Consequently
it does not make immediate sense to define the compactified relative
Picard schemes of $\mathcal{E}\rightarrow\mathbb{P}^5$. However, we
can define $Y^k$ to be a component of the Mukai moduli space of stable
sheaves on $S$; its generic element will be a torsion sheaf obtained
by pushing forward a degree $k$ line bundle on a smooth curve in the
family $\mathcal{E}\rightarrow\mathbb{P}^5$.

We now consider two cases:
\begin{enumerate}
\item when $k=2m+1$ is odd $Y^{2m+1}$ is smooth and compact, and
deformation equivalent to $S^{[5]}$; $Y^{2m+1}$ can therefore be
regarded as a smooth compactification of
$\mathrm{Pic}^{2m+1}(\mathcal{E}/\mathbb{P}^5)$,
\item when $k=2m$ is even $Y^{2m}$ is non-compact; we can compactify
to $\overline{Y}^{2m}$ by adding semi-stable sheaves but this
introduces singularities.
\end{enumerate}
By Hironaka's theorem $\overline{Y}^{2m}$ can be desingularized, but
to find a holomorphic symplectic resolution is more difficult. 

Recently O'Grady~\cite{ogrady99} produced a new ten dimensional
irreducible holomorphic symplectic manifold by showing that a similar
singular moduli space of semi-stable sheaves $\mathcal{M}$ admits a
symplectic desingularization $\widetilde{\mathcal{M}}$. In fact the
singular spaces $\mathcal{M}$ and $\overline{Y}^6$ are birational, and
the structure of their singularities agree, at least
locally~\cite{lehn03}. So we expect that there exists a similar
symplectic desingularization $\widetilde{Y}^6$ of $\overline{Y}^6$,
which would be birational to O'Grady's example
$\widetilde{\mathcal{M}}$ (and hence also deformation equivalent to
it, by a result of Huybrechts~\cite{huybrechts99i}). The space
$\widetilde{Y}^6$ could be regarded as a smooth compactification of 
$\mathrm{Pic}^6(\mathcal{E}/\mathbb{P}^5)$. 

Fix now a specific example from the first case: $Y^5$ for
instance. For a smooth genus five curve in the family
$\mathcal{E}\rightarrow\mathbb{P}^5$, the fibres $\mathrm{Pic}^5$ and
$\mathrm{Pic}^6$ of $Y^5\rightarrow\mathbb{P}^5$ and
$\widetilde{Y}^6\rightarrow\mathbb{P}^5$ are isomorphic smooth abelian
varieties. Denote by $U^5\subset Y^5$ and $U^6\subset\widetilde{Y}^6$
the Zariski open subsets consisting of the union of these smooth
fibres; then in fact $U^5$ and $U^6$ are locally isomorphic
fibrations, but globally we expect a twist.

We summarize this in the following diagram:
$$\begin{array}{ccccccc}
S^{[5]} & \stackrel{\mbox{deform}}{\rightsquigarrow} & Y^5 & &
\widetilde{Y}^6 & \stackrel{\mbox{deform}}{\rightsquigarrow} &
\widetilde{\mathcal{M}} \\
 & & \cup & & \cup & & \\
 & & U^5 & \stackrel{\mbox{twist}}{\rightsquigarrow} & U^6 & & \\
\end{array}$$
The Hilbert scheme $S^{[5]}$ and O'Grady's example
$\widetilde{\mathcal{M}}$ are not deformation equivalent as their
second Betti numbers differ; therefore $Y^5$ and $\widetilde{Y}^6$
cannot be deformation equivalent. It is therefore interesting that
large subsets of them both (i.e.\ $U^5$ and $U^6$) are so closely
related.

The hope of constructing a twisted Fourier-Mukai transform relating
$Y^5$ and $\widetilde{Y}^6$ (and thereby $S^{[5]}$ and
$\widetilde{\mathcal{M}}$) seems remote, as their singular fibres are
quite different: they could contain different numbers of irreducible
components, for instance. Nevertheless, some interesting questions
still arise: can one define a dual fibration of $Y^5$ in a sensible
way? It should agree with $Y^0$ over smooth fibres; how is $Y^0$
related to $Y^6$? Cho, Miyaoka, and Shepherd-Barron~\cite{cmsb02}
described how to construct an abelian fibration $X^0$ with a section
from an abelian fibration $X$ with only stable fibres, such that $X$
is a torsor over $X^0$: over smooth fibres $X^0$ is just
$\mathrm{Pic}^0(\mathrm{Pic}^0(X/B))$, which is then completed by
adding the same singular fibres as occur in $X$. In our case it is
not clear that $Y^5$ does not contain non-reduced, and hence unstable,
fibres (which might occur over non-reduced curves in the family
$\mathcal{E}\rightarrow\mathbb{P}^5$).

\subsection{The generalized Kummer four-fold}

We return to irreducible holomorphic symplectic four-folds. Currently
just two examples are known: the Hilbert scheme of two points on a K3
surface, that we saw in Subsection 5.1, and the generalized Kummer
four-fold $K_4$. Beginning with an abelian surface $A$, the Hilbert
scheme $A^{[3]}$ of three points on $A$ is a minimal resolution of the
symmetric product $\mathrm{Sym}^3A$. It is holomorphic symplectic but
not irreducible. However, there is a composition of maps
$$A^{[3]}\rightarrow\mathrm{Sym}^3A\rightarrow A$$
(the second map is addition of the three points in $A$) and all fibres
of this composition are isomorphic, irreducible holomorphic symplectic
four-folds, which we denote by $K_4$ (see
Beauville~\cite{beauville83}).

\begin{example}
Suppose that $A$ is polarized by a smooth genus four curve $C$ (this
is a type $(1,3)$ polarization). The curve $C$ moves in a
two dimensional linear system $\mathbb{P}^2$. The induced map
$A\rightarrow(\mathbb{P}^2)^{\vee}$ is six-to-one, ramified over a
degree eighteen curve (see Lange and Sernesi~\cite{ls02}), and the
family $\mathcal{C}\rightarrow\mathbb{P}^2$ of curves are pull-back of
lines. If $A$ is chosen to be otherwise generic, then the curves in
this family are all irreducible and we can define the compactified
relative Picard schemes
$Z^k:=\overline{\mathrm{Pic}}^k(\mathcal{C}/\mathbb{P}^2)$, which are
six-folds. 

For each $k$, there is a map $\pi_k:Z^k\rightarrow A$ which is just
the Albanese map of $Z^k$. The kernel $Y^k$ is a smooth four-fold, and
it was shown by Debarre~\cite{debarre99} that $Y^k$ is deformation
equivalent to $K_4$. Restricting $\pi_k$ to a fibre of $Z^k$ gives a
surjective map from a principally polarized four dimensional abelian
variety (or a degeneration thereof) to the $(1,3)$ polarized abelian
surface $A$; the kernel of this map is thus a $(1,3)$ polarized
abelian surface (or degeneration thereof). It follows that $Y^k$ is an
abelian fibration over $\mathbb{P}^2$ whose generic fibre is a
$(1,3)$ polarized abelian surface. 
\end{example}

Unlike the examples we have already seen, whose generic fibres were
principally polarized Jacobians, the fibres of $Y^k$ will not be
self-dual. Assuming that the fibres of $Y^k$ are irreducible, we can
define the compactified relative Picard scheme
$P:=\overline{\mathrm{Pic}}^0(Y^k/\mathbb{P}^2)$. Since this is the
dual fibration, even around smooth fibres it won't be isomorphic to
$Y^k$, only isogeneous. Numerous questions now arise: is $P$ a
deformation of $K_4$? Is it even smooth or holomorphic symplectic?

If $P$ is a deformation of $K_4$, then we observe that the exponential
long exact sequence gives
$$\mathrm{H}^2(P,\mathcal{O})\rightarrow\mathrm{H}^2(P,\mathcal{O}^*)\rightarrow\mathrm{H}^3(P,\mathbb{Z})\rightarrow
0$$
where the first term is $\mathbb{C}$ and the third term is
$\mathbb{Z}^{\oplus 8}$. Thus the space of gerbes
$\mathrm{H}^2(P,\mathcal{O}^*)$ contains an infinite number of
connected components. Can we construct torsors $Y_{\beta}$ for all
gerbes $\beta$? If so, and if $\beta_1$ and $\beta_2$ belong to
different connected components of $\mathrm{H}^2(P,\mathcal{O}^*)$,
then it would appear that $Y_{\beta_1}$ cannot be deformed through
abelian fibrations to $Y_{\beta_2}$. Whether $Y_{\beta_1}$ and
$Y_{\beta_2}$ are related through any deformation at all is another
interesting question.

\bibliographystyle{amsalpha}

\end{document}